\documentclass[12pt]{article}
\usepackage{latexsym,amssymb,amsmath}
\begin{document}
\baselineskip=20pt

\newcommand{\la}{\langle}
\newcommand{\ra}{\rangle}
\newcommand{\psp}{\vspace{0.4cm}}
\newcommand{\pse}{\vspace{0.2cm}}
\newcommand{\ptl}{\partial}
\newcommand{\dlt}{\delta}
\newcommand{\sgm}{\sigma}
\newcommand{\al}{\alpha}
\newcommand{\be}{\beta}
\newcommand{\G}{\Gamma}
\newcommand{\gm}{\gamma}
\newcommand{\vs}{\varsigma}
\newcommand{\Lmd}{\Lambda}
\newcommand{\lmd}{\lambda}
\newcommand{\td}{\tilde}
\newcommand{\vf}{\varphi}
\newcommand{\yt}{Y^{\nu}}
\newcommand{\wt}{\mbox{wt}\:}
\newcommand{\rd}{\mbox{Res}}
\newcommand{\ad}{\mbox{ad}}
\newcommand{\stl}{\stackrel}
\newcommand{\ol}{\overline}
\newcommand{\ul}{\underline}
\newcommand{\es}{\epsilon}
\newcommand{\dmd}{\diamond}
\newcommand{\clt}{\clubsuit}
\newcommand{\vt}{\vartheta}
\newcommand{\ves}{\varepsilon}
\newcommand{\dg}{\dagger}
\newcommand{\tr}{\mbox{Tr}}
\newcommand{\ga}{{\cal G}({\cal A})}
\newcommand{\hga}{\hat{\cal G}({\cal A})}
\newcommand{\Edo}{\mbox{End}\:}
\newcommand{\for}{\mbox{for}}
\newcommand{\kn}{\mbox{ker}}
\newcommand{\Dlt}{\Delta}
\newcommand{\rad}{\mbox{Rad}}
\newcommand{\rta}{\rightarrow}
\newcommand{\mbb}{\mathbb}
\newcommand{\lra}{\Longrightarrow}
\newcommand{\X}{{\cal X}}
\newcommand{\Y}{{\cal Y}}
\newcommand{\Z}{{\cal Z}}
\newcommand{\U}{{\cal U}}
\newcommand{\V}{{\cal V}}
\newcommand{\W}{{\cal W}}
\setlength{\unitlength}{3pt}

\begin{center}{\Large \bf Polynomial Representation of $E_6$ and Its}\end{center}
\begin{center}{\Large \bf  Combinatorial and PDE Implications}\footnote
{2000 Mathematical Subject Classification. Primary 17B10, 17B25;
Secondary 17B01.}
\end{center}
\vspace{0.2cm}

\begin{center}{\large Xiaoping Xu}\end{center}
\begin{center}{Institute of Mathematics, Academy of Mathematics \& System Sciences}\end{center}
\begin{center}{Chinese Academy of Sciences, Beijing 100190, P.R. China}
\footnote{Research supported
 by China NSF 10871193}\end{center}

\vspace{0.6cm}

 \begin{center}{\Large\bf Abstract}\end{center}

\vspace{1cm} {\small  In this paper, we use partial differential
equations to find the decomposition of the polynomial algebra over
the basic irreducible module of $E_6$ into a sum of irreducible
submodules. It turns out that the cubic polynomial invariant
corresponding to the Dicksons' invariant trilinear  form is the
unique fundamental invariant. Moreover, we obtain a combinatorial
identity saying that the dimensions of certain irreducible modules
of $E_6$ are correlated by the binomial coefficients of twenty-six.
Furthermore, we find all the polynomial solutions for the invariant
differential operator corresponding to the Dickson trilinear form in
terms of the irreducible submodules.}

\section{Introduction}

The $E_6$ Lie algebra and group are popular mathematical objects
with broad applications. Dickson [13] (1901) first realized that
there exists an $E_6$-invariant trilinear form on its 27-dimensional
basic irreducible module. The 78-dimensional simple Lie algebra of
type $E_6$ can be realized by all the derivations and multiplication
operators with trace zero on the 27-dimensional exceptional simple
Jordan algebra (e.g., cf. [1], [30]). Aschbacher [5] used the
Dickson form to study the subgroup structure of the $E_6$ group.
Bion-Nadal [8] proved that the $E_6$ Coxeter graph can be realized
as a principal graph of subfactor of the hyperfinite $\Pi_1$ factor.
Brylinski and Kostant [9] obtained a generalized Capelli identity on
the minimal representation of $E_6$. Binegar and Zierau [7] found a
singular representation of $E_6$. Ginzburg [16] proved that the
twisted partial $L$-function on the 27-dimensional representation of
$GE_6(\mbb{C})$ is entire except the points 0 and 1. Iltyakov [21]
showed that the field of invariant rational functions of $E_6$ on
the direct sum of finite copies of the basic module and its dual is
purely transcendental. Suzuki and Wakui [29] studied the
Turaev-Viro-Ocneanu invariant of 3-manifilds derived from the
$E_6$-subfactor. Moreover, Cerchiai and Scotti [11] investigated the
mapping geometry of the $E_6$ group. Furthermore, the $(A_2,G_2)$
duality in $E_6$ was obtained by Rubenthaler [28].

Okamoto and Marshak [27] constructed a grand unification preon model
with $E_6$ metacolor. The $E_6$ Lie algebra was used in [18] to
explain the degeneracies encountered in the genetic code as the
result of a sequence of symmetry breakings that have occurred during
its evolution. Wang [31] identified Geoner's model with the twisted
LG model and $E_6$ singlets. Morrison, Pieruschka and Wybourne [26]
constructed the $E_6$ interacting boson model. Berglund, Candelas et
al. [6] studied instanton contrbutions to the masses and couplings
of $E_6$ singlets. Haba and Matsuoka [17] found large lepton flavor
mixing in the $E_6$-type unification models.
 Ghezelbash, Shafiekhani and Abolbasani [15] derived explicitly a set of Picard-Fuchs
 equations of $N=2$ supersymmetric $E_6$ Yang-Mills
theory. Anderson and Bla$\check z$ek [2-4] found certain
Clebsch-Gordan coefficients in connection with the $E_6$ unification
model building. Fern\'{a}ndez-N\'{u}$\td{a}$ez, Garcia-Fuertes and
Perelomov [14] used the quantum Calogero-Sutherland model
corresponding to the root system of $E_6$ to calculate
Clebach-Gordan series for this algebra. Howl and King [20] proposed
a minimal $E_6$ supersymmetric standard model which allows Planck
scale unification, provides a solution to the $\mu$ problem and
predicts a new $Z'$. Das and Laperashvili [12] studied Preon model
related to family replicated $E_6$ unification.

The motivation of this work is to understand the functional impact
of the $E_6$ Lie algebra (group). We  use partial differential
equations to decompose the polynomial algebra over the basic
irreducible module of $E_6$ into a sum of irreducible submodules.
Consequently, the cubic polynomial invariant corresponding to the
Dickson's invariant trilinear  form is the unique fundamental
invariant (any other invariant is a polynomial in it). Moreover, we
obtain a combinatorial identity, which  says that the dimensions of
certain irreducible modules of $E_6$ are correlated by the binomial
coefficients of twenty-six. Furthermore, we find all the polynomial
solutions for the invariant differential operator corresponding to
the Dickson trilinear form in terms of the irreducible submodules.
Below we give a more detailed introduction to our results.

It has been known for may years that the representation theory of
Lie algebra is closely related to combinatorial identities.
Macdonald [25] generalized the Weyl denominator identities for
finite root systems to those for infinite affine root systems, which
are now known as the Macdonald's identities. We present here a
consequence  of the Macdonald's identities taken from Kostant's work
[22]. Let ${\cal G}$ be a finite-dimensional simple Lie algebra over
the field $\mbb{C}$ of complex numbers. Denote by $\Lmd^+$ the set
of dominant weights of ${\cal G}$ and by $V(\lmd)$ the
finite-dimensional irreducible ${\cal G}$-module with highest weight
$\lmd$. It is known that the Casimir operator takes a constant
$c(\lmd)$ on $V(\lmd)$. Macdonald's Theorem implies that there
exists a map $\chi:\Lmd^+\rta\{-1,0,1\}$ such that the following
identity holds:
$$(\prod_{n=1}^\infty(1-q^n))^{\dim {\cal
G}}=\sum_{\lmd\in\Lmd^+}\chi(\lmd)(\dim
V(\lmd))q^{c(\lmd)}.\eqno(1.1)$$ Kostant [23]  connected the above
identity to the abelian subalgebras of ${\cal G}$.

Denote by $E_{r,s}$ the square matrix with 1 as its $(r,s)$-entry
and 0 as the others. The orthogonal Lie algebra
$$o(n,\mbb{C})=\sum_{1\leq r<s\leq
n}\mbb{C}(E_{r,s}-E_{s,r}).\eqno(1.2)$$ It acts on the polynomial
algebra ${\cal A}=\mbb{C}[x_1,...,x_n]$ by
$$(E_{r,s}-E_{s,r})|_{\cal
A}=x_r\ptl_{x_s}-x_s\ptl_{x_r}.\eqno(1.3)$$ Denote by ${\cal A}_k$
the subspace of homogeneous polynomials in ${\cal A}$ with degree
$k$. When $n\geq 3$, it is well known that  the subspace of harmonic
polynomials
$${\cal H}_k=\{f\in{\cal A}_k\mid
(\ptl_{x_1}^2+\cdots+\ptl_{x_n}^2)(f)=0\}\eqno(1.4)$$ forms an
irreducible $o(n,\mbb{C})$-module and
$${\cal A}_k={\cal H}_k\oplus(x_1^2+x_2^2+\cdots+x_n^2){\cal
A}_{k-2}.\eqno(1.5)$$

Recall the special linear Lie algebra
$$sl(n,\mbb{C})=\sum_{r\neq
s}\mbb{C}E_{r,s}+\sum_{r=1}^{n-1}\mbb{C}(E_{r,r}-E_{r+1,r+1}).\eqno(1.6)$$
Let ${\cal A} $ be the polynomial algebra in
$x_1,...,x_n,y_1,...,y_n.$  Define a representation of
$sl(n,\mbb{C})$ on ${\cal A}$ via
$$E_{r,s}|_{\cal
A}=x_r\ptl_{x_s}-y_s\ptl_{y_r}.\eqno(1.7)$$ Denote by $\mbb{N}$ the
additive semigroup of nonnegative integers. Define
$$x^\al=x_1^{\al_1}x_2^{\al_2}\cdots
x_n^{\al_n}\qquad\for\;\;\al=(\al_1,...,\al_n)\in\mbb{N}^{\:n}.\eqno(1.8)$$
 Set
$${\cal
A}_{\ell_1,\ell_2}=\sum_{\al,\be\in\mbb{N}^{\:n};\;|\al|=\ell_1,\;|\be|=\ell_2}
\mbb{R}x^\al y^\be\qquad\for\;\;\ell_1,\ell_2\in\mbb{N}\eqno(1.9)$$
and define
$${\cal H}_{\ell_1,\ell_2}=\{u\in {\cal
A}_{\ell_1,\ell_2}\mid
u_{x_1y_1}+u_{x_2y_2}+\cdots+u_{x_ny_n}=0\}.\eqno(1.10)$$  We proved
in [32] that ${\cal H}_{\ell_1,\ell_2}$ forms an irreducible
$sl(n,\mbb{C})$-submodule and
$${\cal A}_{\ell_1,\ell_2}={\cal H}_{\ell_1,\ell_2}\oplus
(x_1y_1+x_2y_2+\cdots+x_ny_n){\cal
A}_{\ell_1-1,\ell_2-1},\eqno(1.11)$$ where
$x_1y_1+x_2y_2+\cdots+x_ny_n$ is an  $sl(n,\mbb{C})$-invariant. More
importantly, an explicit basis for each ${\cal H}_{\ell_1,\ell_2}$
was constructed in [32].

Now we assume that ${\cal A}$ is the polynomial algebra in
$x_1,...,x_7$. There exists an action of the simple Lie algebra
${\cal G}^{G_2}$ of type $G_2$ on ${\cal A}$, which keeps
$x_1^2+x_2x_5+x_3x_6+x_4x_7$ invariant (e.g., cf. [32]). Again we
denote by ${\cal A}_k$ the subspace of polynomials of degree $k$ in
${\cal A}$ and define
$$\tilde{\cal H}_k=\{u\in{\cal A}_k\mid
u_{x_1x_1}+u_{x_2x_5}+u_{x_3x_6}+u_{x_4x_7}=0\}.\eqno(1.12)$$ We
showed in [32]  that $\tilde{\cal H}_k$ forms an irreducible ${\cal
G}^{G_2}$-submodule and
$${\cal A}_k=\tilde{\cal H}_k\oplus(x_1^2+x_2x_5+x_3x_6+x_4x_7){\cal
A}_{k-2}.\eqno(1.13)$$ Moreover, an explicit basis for each
$\tilde{\cal H}_k$ was constructed in [32]. Furthermore, Luo [24]
generalized the result (1.5) to certain noncanonical polynomial
representations of $o(n,\mbb{C})$ and found a basis for the
irreducible submodules by using our methods in [32].  The methods
also enabled Cao [10] to obtain an explicit irreducible
decomposition with respect to the real orthogonal Lie algebra for
the space of polynomial solutions of the Navier equations in
elasticity.

The Dynkin diagram of $E_6$ is as follows:

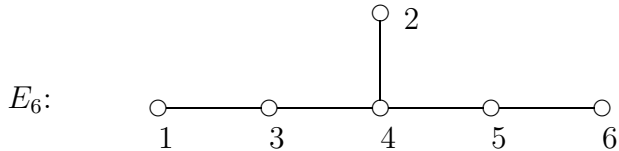
\begin{picture}(80,20)
\put(2,0){$E_6$:}\put(21,0){\circle{2}}\put(21,
-5){1}\put(22,0){\line(1,0){12}}\put(35,0){\circle{2}}\put(35,
-5){3}\put(36,0){\line(1,0){12}}\put(49,0){\circle{2}}\put(49,
-5){4}\put(49,1){\line(0,1){10}}\put(49,12){\circle{2}}\put(52,10){2}\put(50,0){\line(1,0){12}}
\put(63,0){\circle{2}}\put(63,-5){5}\put(64,0){\line(1,0){12}}\put(77,0){\circle{2}}\put(77,
-5){6}
\end{picture}
\vspace{0.7cm}

\noindent Denote by $\lmd_i$ the $i$th fundamental weight of $E_6$
with respect to the above labeling. Let $V$ be the 27-dimensional
irreducible $E_6$-module of highest weight $\lmd_1$. Denote by
${\cal A}$ the polynomial algebra (equivalently, symmetric tensor)
over $V$ and  by ${\cal A}_m$ the subspace of homogeneous polynomial
with degree $m $. A {\it singular vector} in ${\cal A}$ is a weight
vector annihilated by positive root vectors. The following is the
main theorem of this paper:\psp

{\bf Main Theorem}: {\it Any singular vector is a monomial in a
linear singular vector $x_1$ of weight $\lmd_1$, a quadratic
singular vector $\zeta_1$ of weight $\lmd_6$ and a cubic singular
vector $\eta$ of weight $0$. In particular, $\eta$ is the unique
fundamental invariant and the following combinatorial identity
holds:
$$(1-q)^{26}\sum_{m_1,m_2=0}^\infty (\mbox{\it dim}\:
V(m_1\lmd_1+m_2\lmd_6))q^{m_1+2m_2}=1+q+q^2.\eqno(1.14)$$ Let ${\cal
D}$ be the unique constant-coefficient fundamental invariant
differential operator dual to $\eta$. Denote by $L(m_1,m_2)$ the
$E_6$-submodule generated by $x_1^{m_1}\zeta_1^{m_2}$. Then
$$\Phi_m=\{f\in{\cal A}_m\mid {\cal
D}(f)=0\}=\bigoplus_{i=0}^{[\!|m/2|\!]}L(m-2i,i)\eqno(1.15)$$ and
$${\cal A}_m=\Phi_m\oplus \eta {\cal A}_{m-3}.\eqno(1.16)$$}\pse

Note that the identities (1.1) and (1.14) are dimensional properties
of irreducible submodules. Our identity (1.14) says that the
dimensions of the irreducible $E_6$-modules $V(m_1\lmd_1+m_2\lmd_6)$
are correlated by the binomial coefficients of twenty-six. Moreover,
$\eta$ is the cubic polynomial invariant corresponding to the
Dickson's invariant trilinear form (cf. [13]) and ${\cal D}$ is the
invariant differential operator corresponding to the Dickson form.
The equation (1.16) is exactly a cubic generalization of the
quadratic ones in (1.5), (1.11) and (1.13). The fundamental
difference is that our subspace $\Phi_m$ of homogeneous polynomial
solutions is a sum of $[\!|m/2|\!]+1$ irreducible submodules.

In Section 2, we explicitly construct the 27-dimensional basic
representation of $E_6$ in terms of differential operators  via the
root lattice construction of the $E_7$ simple Lie algebra. The proof
of the main theorem is given in Section 3.

 \section{Basic Representation of $E_6$}

In this section, we will explicitly construct the 27-dimensional
basic irreducible representation of $E_6$.

 For convenience, we will use the notion
$$\ol{i,i+j}=\{i,i+1,i+2,...,i+j\}\eqno(2.1)$$
for integer $i$ and positive integer $j$ throughout this paper. We
start with the root lattice construction of the simple Lie algebra
of type $E_7$. As we all known, the Dynkin diagram of $E_7$ is as
follows:

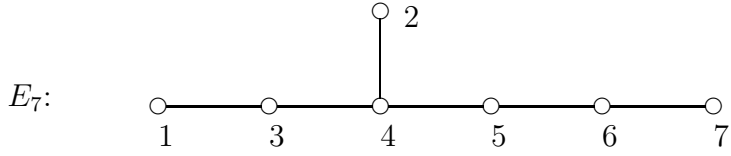
\begin{picture}(93,20)
\put(2,0){$E_7$:}\put(21,0){\circle{2}}\put(21,
-5){1}\put(22,0){\line(1,0){12}}\put(35,0){\circle{2}}\put(35,
-5){3}\put(36,0){\line(1,0){12}}\put(49,0){\circle{2}}\put(49,
-5){4}\put(49,1){\line(0,1){10}}\put(49,12){\circle{2}}\put(52,10){2}\put(50,0){\line(1,0){12}}
\put(63,0){\circle{2}}\put(63,-5){5}\put(64,0){\line(1,0){12}}\put(77,0){\circle{2}}\put(77,
-5){6}\put(78,0){\line(1,0){12}}\put(91,0){\circle{2}}\put(91,
-5){7}
\end{picture}
\vspace{0.7cm}

 \noindent Let $\{\al_i\mid i\in\ol{1,7}\}$ be the
simple positive roots corresponding to the vertices in the diagram,
and let $\Phi_{E_7}$ be the root system of $E_7$. Set
$$Q_{E_7}=\sum_{i=1}^7\mbb{Z}\al_i,\eqno(2.2)$$ the root lattice of type
$E_7$. Denote by $(\cdot,\cdot)$ the symmetric $\mbb{Z}$-bilinear
form on $Q_{E_7}$ such that
$$\Phi_{E_7}=\{\al\in Q_{E_7}\mid (\al,\al)=2\}.\eqno(2.3)$$
Define $F(\cdot,\cdot):\; Q_{E_7}\times  Q_{E_7}\rta \{\pm 1\}$ by
$$F(\sum_{i=1}^7k_i\al_i,\sum_{j=1}^7l_j\al_j)=(-1)^{\sum_{i=1}^7k_il_i+\sum_{7\geq i>j\geq 1}k_il_j
(\al_i,\al_j)},\qquad k_i,l_j\in\mbb{Z}.\eqno(2.4)$$ Then for
$\al,\be,\gm\in  Q_{E_7}$,
$$F(\al+\be,\gm)=F(\al,\gm)F(\be,\gm),\;\;F(\al,\be+\gm)=F(\al,\be)F(\al,\gm),\eqno(2.5)
$$
$$F(\al,\be)F(\be,\al)^{-1}=(-1)^{(\al,\be)},\;\;F(\al,\al)=(-1)^{(\al,\al)/2}.
\eqno(2.6)$$ In particular,
$$F(\al,\be)=-F(\be,\al)\qquad
\mbox{if}\;\;\al,\be,\al+\be\in \Phi_{E_7}.\eqno(2.7)$$

Denote
$$H_{E_7}=\sum_{i=1}^7\mbb{C}\al_i.\eqno(2.8)$$
The simple Lie algebra of type $E_7$ is
$${\cal
G}^{E_7}=H_{E_7}\oplus\bigoplus_{\al\in
\Phi_{E_7}}\mbb{C}E_{\al}\eqno(2.9)$$ with the Lie bracket
$[\cdot,\cdot]$ determined by:
 $$[H_{E_7},H_{E_7}]=0,\;\;[h,E_{\al}]=(h,\al)E_{\al},\;\;[E_{\al},E_{-\al}]=-\al,
 \eqno(2.10)$$
 $$[E_{\al},E_{\be}]=\left\{\begin{array}{ll}0&\mbox{if}\;\al+\be\not\in \Phi_{E_7},\\
 F(\al,\be)E_{\al+\be}&\mbox{if}\;\al+\be\in\Phi_{E_7}.\end{array}\right.\eqno(2.11)$$
for $\al,\be\in\Phi_{E_7}$ and $h\in H_{E_7}$.

Note that
$$Q_{E_6}=\sum_{i=1}^6\mbb{Z}\al_i\subset Q_{E_7}\eqno(2.12)$$
is the root lattice of $E_6$ and
$$\Phi_{E_6}=Q_{E_6}\bigcap \Phi_{E_7}\eqno(2.13)$$
is the root system of $E_6$. Set
$$H_{E_6}=\sum_{i=1}^6\mbb{C}\al_i.\eqno(2.14)$$
Then the subalgebra
$${\cal G}^{E_6}=H_{E_6}\oplus\bigoplus_{\al\in
\Phi_{E_6}}\mbb{C}E_{\al}\eqno(2.15)$$ of ${\cal G}^{E_7}$ is
exactly the simple Lie algebra of ${\cal G}^{E_6}$. Denote by
$\Phi_{E_6}^+$ the set of positive roots of $E_6$ and by
$\Phi_{E_7}^+$ the set of positive roots of $E_7$. The elements of
$\Phi_{E_6}^+$ are:
$$\al_1+2\al_2+2\al_3+3\al_4+2\al_5+\al_6,\eqno(2.16)$$
$$\{\al_1+\sum_{r=3}^j\al_r\mid j\in\ol{2,6}\}\bigcup
\{\sum_{r=i+1}^j\al_r\mid 2\leq i<j\leq 6\},\eqno(2.17)$$
$$\{\sum_{s=2}^j\al_s+\sum_{t=4}^k\al_t\mid 2\leq j<k\leq
6\}\eqno(2.18)$$ and
$$\{\sum_{\iota=1}^i\al_\iota+
\sum_{s=3}^j\al_s+\sum_{t=4}^k\al_t\mid 2\leq i< j<k\leq
6\}.\eqno(2.19)$$

Denote by $\bar\Phi_{E_7}^+$ the set of the following positive
roots:
$$\al_1+\sum_{r=3}^7\al_r,\qquad\al_3+2\al_4+\al_5+\sum_{i=1}^6\al_i+\sum_{r=1}^7\al_r,\eqno(2.20)$$
$$\{2\sum_{s=1}^6\al_s-\al_1+\al_4-\al_6+\sum_{r=i+1}^7\al_r\mid
i\in\ol{1,6}\} ,\qquad \{\sum_{r=i+1}^7\al_r\mid
i\in\ol{2,6}\},\eqno(2.21)$$
$$\{\sum_{s=2}^j\al_s+\sum_{t=4}^7\al_t\mid j\in\ol{2,6}\},\qquad \{\sum_{\iota=1}^i\al_\iota+
\sum_{s=3}^j\al_s+\sum_{t=4}^7\al_t\mid 2\leq i< j\leq
6\}.\eqno(2.22)$$ Then
$$\Phi_{E_7}^+=\Phi_{E_6}^+\bigcup \bar\Phi_{E_7}^+.\eqno(2.23)$$

In particular,
$$V=\sum_{\be\in \bar\Phi_{E_7}^+}\mbb{C}E_\be\eqno(2.24)$$
forms the 27-dimensional basic ${\cal G}^{E_6}$-module of highest
weight $\lmd_1$ with the representation $\mbox{ad}_{{\cal
G}^{E_7}}$. Denote
$$x_1=E_{\al_3+2\al_4+\al_5+\sum_{i=1}^6\al_i+\sum_{r=1}^7\al_r},\qquad
x_2=E_{2\sum_{s=1}^6\al_s-\al_1+\al_4-\al_6+\sum_{r=3}^7\al_r},\eqno(2.25)$$
$$x_3=E_{2\sum_{s=1}^6\al_s-\al_1+\al_4-\al_6+\sum_{r=4}^7\al_r},\qquad
x_4=E_{2\sum_{s=1}^6\al_s-\al_1-\al_6+\sum_{r=4}^7\al_r}\eqno(2.26)$$
$$x_5=E_{2\sum_{s=1}^6\al_s-\al_1+\al_4+\al_7},\qquad x_6
=E_{\sum_{\iota=1}^5\al_\iota+\sum_{s=3}^6\al_s+\sum_{t=4}^7\al_t},\eqno(2.27)$$
$$x_7 =E_{\al_4+\sum_{s=3}^6\al_s+\sum_{t=1}^7\al_t},\qquad
x_8=E_{2\sum_{s=1}^6\al_s-\al_1+\al_4-\al_6+\al_7},\eqno(2.28)$$
$$x_{9}=E_{\sum_{s=3}^6\al_s+\sum_{t=1}^7\al_t},\qquad
x_{10}=E_{\al_4+\sum_{s=3}^5\al_s+\sum_{t=1}^7\al_t},\eqno(2.29)$$
$$x_{11}
=E_{\sum_{s=4}^6\al_s+\sum_{t=1}^7\al_t},\qquad x_{12}
=E_{\sum_{s=3}^5\al_s+\sum_{t=1}^7\al_t},\qquad
 x_{13}=E_{\al_3+\al_4+\sum_{t=1}^7\al_t},\eqno(2.30)
$$
$$x_{14}=E_{\sum_{r=2}^6 \al_r+\sum_{s=4}^7\al_s}, \qquad x_{15}
=E_{\al_4+\al_5+\sum_{t=1}^7\al_t},\qquad x_{16}
=E_{\al_4+\sum_{t=1}^7\al_t} ,\eqno(2.31)$$ $$ x_{17}
=E_{\al_4+\al_5+\sum_{t=2}^7\al_t},\qquad
x_{18}=E_{\sum_{t=1}^7\al_t} ,\qquad x_{19}
=E_{\al_4+\sum_{t=2}^7\al_t},\eqno(2.32)$$
$$x_{20}=E_{\al_1+\sum_{t=3}^7\al_t},\qquad x_{21}=E_{\sum_{t=2}^7\al_t},\qquad
x_{22}=E_{\al_2+\sum_{t=4}^7\al_t},\qquad
x_{23}=E_{\sum_{t=3}^7\al_t},\eqno(2.33)$$
$$x_{24}=E_{\sum_{t=4}^7\al_t},\qquad
x_{25}=E_{\sum_{t=5}^7\al_t},\qquad x_{26}=E_{\al_6+\al_7},\qquad
x_{27}=E_{\al_7}.\eqno(2.34)$$ Then $\{x_i\mid i\in\ol{1,27}\}$
forms a basis of $V$.

Under the above basis
$$E_{\al_1}|_V=-x_1\ptl_{x_2}+x_{11}\ptl_{x_{14}}+x_{15}\ptl_{x_{17}}
+x_{16}\ptl_{x_{19}}+x_{18}\ptl_{x_{21}}+x_{20}\ptl_{x_{23}},\eqno(2.35)$$
$$E_{\al_2}|_V=-x_4\ptl_{x_6}-x_5\ptl_{x_7}-x_8\ptl_{x_{10}}
+x_{18}\ptl_{x_{20}}+x_{21}\ptl_{x_{23}}+x_{22}\ptl_{x_{24}},\eqno(2.36)$$
$$E_{\al_3}|_V=-x_2\ptl_{x_3}+x_9\ptl_{x_{11}}
+x_{12}\ptl_{x_{15}}+x_{13}\ptl_{x_{16}}+x_{21}\ptl_{x_{22}}+x_{23}\ptl_{x_{24}},\eqno(2.37)$$
$$E_{\al_4}|_V=-x_3\ptl_{x_4}-x_7\ptl_{x_9}-x_{10}\ptl_{x_{12}}
-x_{16}\ptl_{x_{18}}-x_{19}\ptl_{x_{21}}+x_{24}\ptl_{x_{25}},\eqno(2.38)$$
$$E_{\al_5}|_V=-x_4\ptl_{x_5}-x_6\ptl_{x_7}-x_{12}\ptl_{x_{13}}-x_{15}\ptl_{x_{16}}
-x_{17}\ptl_{x_{19}}+x_{25}\ptl_{x_{26}},\eqno(2.39)$$
$$E_{\al_6}|_V=-x_5\ptl_{x_8}-x_7\ptl_{x_{10}}-x_9\ptl_{x_{12}}-x_{11}
\ptl_{x_{15}} -x_{14}\ptl_{x_{17}}+x_{26}\ptl_{x_{27}},\eqno(2.40)$$
$$E_{\al_1+\al_3}|_V=x_1\ptl_{x_3}-x_9\ptl_{x_{14}}-x_{12}\ptl_{x_{17}}
-x_{13}\ptl_{x_{19}}+x_{18}\ptl_{x_{22}}+x_{20}\ptl_{x_{24}},\eqno(2.41)$$
$$E_{\al_2+\al_4}|_V=-x_3\ptl_{x_6}+x_5\ptl_{x_9}+x_8\ptl_{x_{12}}+x_{16}\ptl_{x_{20}}
+x_{19}\ptl_{x_{23}}+x_{22}\ptl_{x_{25}},\eqno(2.42)$$
$$E_{\al_3+\al_4}|_V=x_2\ptl_{x_4}+x_7\ptl_{x_{11}}+x_{10}\ptl_{x_{15}}-x_{13}\ptl_{x_{18}}
+x_{19}\ptl_{x_{22}}+x_{23}\ptl_{x_{25}},\eqno(2.43)$$
$$E_{\al_4+\al_5}|_V=x_3\ptl_{x_5}-x_6\ptl_{x_9}+x_{10}\ptl_{x_{13}}-x_{15}\ptl_{x_{18}}
-x_{17}\ptl_{x_{21}}+x_{24}\ptl_{x_{26}},\eqno(2.44)$$
$$E_{\al_5+\al_6}|_V=x_4\ptl_{x_8}+x_6\ptl_{x_{10}}-x_9\ptl_{x_{13}}-x_{11}\ptl_{x_{16}}-x_{14}\ptl_{x_{19}}
+x_{25}\ptl_{x_{27}},\eqno(2.45)$$
$$E_{\al_1+\al_3+\al_4}|_V=-x_1\ptl_{x_4}-x_7\ptl_{x_{14}}-x_{10}\ptl_{x_{17}}+x_{13}\ptl_{x_{21}}+x_{16}\ptl_{x_{22}}+x_{20}\ptl_{x_{25}},
\eqno(2.46)$$
$$E_{\al_2+\al_3+\al_4}|_V=x_2\ptl_{x_6}-x_5\ptl_{x_{11}}-x_8\ptl_{x_{15}}+x_{13}\ptl_{x_{20}}
-x_{19}\ptl_{x_{x_{24}}}+x_{21}\ptl_{x_{25}},\eqno(2.47)$$
$$E_{\al_2+\al_4+\al_5}|_V=x_3\ptl_{x_7}+
x_4\ptl_{x_9}-x_8\ptl_{x_{13}}+x_{15}\ptl_{x_{20}}+x_{17}\ptl_{x_{23}}+x_{22}\ptl_{x_{26}},\eqno(2.48)$$
$$E_{\al_3+\al_4+\al_5}|_V=-x_2\ptl_{x_5}+x_6\ptl_{x_{11}}-x_{10}\ptl_{x_{16}}-x_{12}\ptl_{x_{18}}
+x_{17}\ptl_{x_{22}}+x_{23}\ptl_{x_{26}},\eqno(2.49)$$
$$E_{\al_4+\al_5+\al_6}|_V=-x_3\ptl_{x_8}+x_6\ptl_{x_{12}}+x_7\ptl_{x_{13}}-x_{11}\ptl_{x_{18}}
-x_{14}\ptl_{x_{21}} +x_{24}\ptl_{x_{27}},\eqno(2.50)$$
$$E_{\sum_{i=1}^4\al_i}|_V=-x_1\ptl_{x_6}+x_5\ptl_{x_{14}}
+x_8\ptl_{x_{17}}-x_{13}\ptl_{x_{23}}-x_{16}\ptl_{x_{24}}+x_{18}\ptl_{x_{25}},\eqno(2.51)$$
$$E_{\al_1+\sum_{i=3}^5\al_i}|_V=x_1\ptl_{x_5}-x_6\ptl_{x_{14}}+x_{10}\ptl_{x_{19}}+x_{12}\ptl_{x_{21}}+x_{15}\ptl_{x_{22}}
+x_{20}\ptl_{x_{26}},\eqno(2.52)$$
$$E_{\sum_{i=2}^5\al_i}|_V=-x_2\ptl_{x_7}-x_4\ptl_{x_{11}}+x_8\ptl_{x_{16}}+x_{12}\ptl_{x_{20}}
-x_{17}\ptl_{x_{24}}+x_{21}\ptl_{x_{26}},\eqno(2.53)$$
$$E_{\al_2+\sum_{i=4}^6\al_i}|_V=-x_3\ptl_{x_{10}}-x_4\ptl_{x_{12}}-x_5\ptl_{x_{13}}+x_{11}\ptl_{x_{20}}+x_{14}\ptl_{x_{23}}
+x_{22}\ptl_{x_{27}},\eqno(2.54)$$
$$E_{\sum_{i=3}^6\al_i}|_V=x_2\ptl_{x_8}-x_6\ptl_{x_{15}}-x_7\ptl_{x_{16}}-x_9\ptl_{x_{18}}+x_{14}\ptl_{x_{22}}+x_{23}\ptl_{x_{27}},
\eqno(2.55)$$
$$E_{\sum_{i=1}^5\al_i}|_V=x_1\ptl_{x_7}+x_4\ptl_{x_{14}}-x_8\ptl_{x_{19}}-x_{12}\ptl_{x_{23}}
-x_{15}\ptl_{x_{24}}+x_{18}\ptl_{x_{26}},\eqno(2.56)$$
$$E_{\al_1+\sum_{i=3}^6\al_i}|_V=-x_1\ptl_{x_8}+x_6\ptl_{x_{17}}+x_7\ptl_{x_{19}}+x_9\ptl_{x_{21}}
+x_{11}\ptl_{x_{22}}+x_{20}\ptl_{x_{27}}, \eqno(2.57)$$
$$E_{\al_4+\sum_{i=2}^5\al_i}|_V=x_2\ptl_{x_9}-x_3\ptl_{x_{11}}-x_8\ptl_{x_{18}}+x_{10}\ptl_{x_{20}}
-x_{17}\ptl_{x_{25}}+x_{19}\ptl_{x_{26}},\eqno(2.58)$$
$$E_{\sum_{i=2}^6\al_i}|_V=x_2\ptl_{x_{10}}+x_4\ptl_{x_{15}}+x_5\ptl_{x_{16}}+x_9\ptl_{x_{20}}-x_{14}\ptl_{x_{24}}+x_{21}\ptl_{x_{27}},
\eqno(2.59)$$
$$E_{\al_4+\sum_{i=1}^5\al_i}|_V=-x_1\ptl_{x_9}+x_3\ptl_{x_{14}}+x_8\ptl_{x_{21}}-x_{10}\ptl_{x_{23}}
-x_{15}\ptl_{x_{25}}+x_{16}\ptl_{x_{26}},\eqno(2.60)$$
$$E_{\sum_{i=1}^6\al_i}|_V=-x_1\ptl_{x_{10}}
-x_4\ptl_{x_{17}}-x_5\ptl_{x_{19}}-x_9\ptl_{x_{23}}-x_{11}\ptl_{x_{24}}+x_{18}\ptl_{x_{27}},
\eqno(2.61)$$
$$E_{\al_4+\sum_{i=2}^6\al_i}|_V=-x_2\ptl_{x_{12}}+x_3\ptl_{x_{15}}
-x_5\ptl_{x_{18}}+x_7\ptl_{x_{20}}-x_{14}\ptl_{x_{25}}+x_{19}\ptl_{x_{27}},
\eqno(2.62)$$
$$E_{\al_3+\al_4+\sum_{i=1}^5\al_i}|_V=x_1\ptl_{x_{11}}
-x_2\ptl_{x_{14}}-x_8\ptl_{x_{22}}+x_{10}\ptl_{x_{24}}
-x_{12}\ptl_{x_{25}}+x_{13}\ptl_{x_{26}},\eqno(2.63)$$
$$E_{\al_4+\sum_{i=1}^6\al_i}|_V=x_1\ptl_{x_{12}}
-x_3\ptl_{x_{17}}+x_5\ptl_{x_{21}}-x_7\ptl_{x_{23}}-x_{11}\ptl_{x_{25}}+x_{16}\ptl_{x_{27}},
\eqno(2.64)$$
$$E_{\al_4+\al_5+\sum_{i=2}^6\al_i}|_V=x_2\ptl_{x_{13}}-x_3\ptl_{x_{16}}
-x_4\ptl_{x_{18}}+x_6\ptl_{x_{20}}-x_{14}\ptl_{x_{26}}+x_{17}\ptl_{x_{27}},
\eqno(2.65)$$
$$E_{\al_3+\al_4+\sum_{i=1}^6\al_i}|_V=-x_1\ptl_{x_{15}}
+x_2\ptl_{x_{17}}-x_5\ptl_{x_{22}}+x_7\ptl_{x_{24}}-x_9\ptl_{x_{25}}+x_{13}\ptl_{x_{27}},
\eqno(2.66)$$
$$E_{\al_4+\al_5+\sum_{i=1}^6\al_i}|_V=-x_1\ptl_{x_{13}}
+x_3\ptl_{x_{19}}+x_4\ptl_{x_{21}}-x_6\ptl_{x_{23}}-x_{11}\ptl_{x_{26}}+x_{15}\ptl_{x_{27}},
\eqno(2.67)$$
$$E_{\sum_{r=3}^5\al_r+\sum_{i=1}^6\al_i}|_V=x_1\ptl_{x_{16}}
-x_2\ptl_{x_{19}}-x_4\ptl_{x_{22}}+x_6\ptl_{x_{24}}-x_9\ptl_{x_{26}}+x_{12}\ptl_{x_{27}},
\eqno(2.68)$$
$$E_{\al_4+\sum_{r=3}^5\al_r+\sum_{i=1}^6\al_i}|_V=x_1\ptl_{x_{18}}
-x_2\ptl_{x_{21}}+x_3\ptl_{x_{22}}-x_6\ptl_{x_{25}}+x_7\ptl_{x_{26}}-x_{10}\ptl_{x_{27}},
\eqno(2.69)$$
$$E_{\al_2+\al_4+\sum_{r=3}^5\al_r+\sum_{i=1}^6\al_i}|_V=x_1\ptl_{x_{20}}
-x_2\ptl_{x_{23}}+x_3\ptl_{x_{24}}-x_4\ptl_{x_{25}}-x_5\ptl_{x_{26}}-x_8\ptl_{x_{27}}.
\eqno(2.70)$$

Recall that we also view $\al_i$ as the elements of ${\cal G}^{E_7}$
(cf. (2.8) and (2.9)). We write
$$[\al_j,x_i]=a_{i,j}x_i\qquad
\for\;\;i\in\ol{1,27},\;j\in\ol{1,6}.\eqno(2.71)$$ Then the weight
of $x_i$ is $\sum_{j=1}^6a_{i,j}\lmd_j$, where $\lmd_j$ is the $j$th
fundamental weight of ${\cal G}^{E_6}$. We calculate the following
table:
\begin{center}{\bf \large Table 1}\end{center}
\begin{center}\begin{tabular}{|r||r|r|r|r|r|r||r||r|r|r|r|r|r|}\hline
$i$&$a_{i,1}$&$a_{i,2}$&$a_{i,3}$&$a_{i,4}$&$a_{i,5}$&$a_{i,6}$&$i$&$a_{i,1}$&$a_{i,2}$&$a_{i,3}$&$a_{i,4}$&$a_{i,5}$&
$a_{i,6}$
\\\hline\hline 1&1&0&0&0&0&0&2&$-1$&0&1&0&0&0\\\hline 3&0&0&$-1$&1&0&0&4&$0$&1&0&$-1$&1&0
\\\hline 5&0&1&$0$&0&$-1$&1&6&$0$&$-1$&0&$0$&1&0 \\\hline
7&0&$-1$&$0$&1&$-1$&1&8&$0$&1&0&$0$&0&$-1$\\\hline
9&0&$0$&$1$&$-1$&0&1&10&$0$&$-1$&0&$1$&0&$-1$\\\hline
11&1&$0$&$-1$&$0$&0&1&12&$0$&0&1&$-1$&1&$-1$\\\hline
13&0&0&1&0&$-1$&0&14&$-1$&0&0&0&0&1\\\hline
15&1&0&$-1$&0&1&$-1$&16&1&0&$-1$&1&$-1$&0\\\hline
17&$-1$&0&0&0&1&$-1$&18&1&1&$0$&$-1$&0&0\\\hline
19&$-1$&0&0&1&$-1$&0&20&1&$-1$&$0$&0&0&0\\\hline
21&$-1$&1&1&$-1$&0&0&22&0&1&$-1$&0&0&0\\\hline
23&$-1$&$-1$&1&0&0&0&24&0&$-1$&$-1$&1&0&0\\\hline
25&0&0&0&$-1$&1&0&26&0&0&0&$0$&$-1$&1\\\hline 27&0&0&0&0&0&$-1$
&&&&&&&\\\hline\end{tabular}\end{center} In particular,
$$\al_j|_V=\sum_{i=1}^{27}a_{i,j}x_i\ptl_{x_i}\qquad\for\;\;j\in\ol{1,6}.
\eqno(2.72)$$

Finally we have the representation formulas of the negative root
vectors:
$$E_{-\al_1}|_V=x_2\ptl_{x_1}-x_{14}\ptl_{x_{11}}-x_{17}\ptl_{x_{15}}
-x_{19}\ptl_{x_{16}}-x_{21}\ptl_{x_{18}}-x_{23}\ptl_{x_{20}},\eqno(2.73)$$
$$E_{-\al_2}|_V=x_6\ptl_{x_4}+x_7\ptl_{x_5}+x_{10}\ptl_{x_8}
-x_{20}\ptl_{x_{18}}-x_{23}\ptl_{x_{21}}-x_{24}\ptl_{x_{22}},\eqno(2.74)$$
$$E_{-\al_3}|_V=x_3\ptl_{x_2}-x_{11}\ptl_{x_9}
-x_{15}\ptl_{x_{12}}-x_{16}\ptl_{x_{13}}-x_{22}\ptl_{x_{21}}
-x_{24}\ptl_{x_{23}},\eqno(2.75)$$
$$E_{-\al_4}|_V=x_4\ptl_{x_3}+x_9\ptl_{x_7}+x_{12}\ptl_{x_{10}}
+x_{18}\ptl_{x_{16}}+x_{21}\ptl_{x_{19}}-x_{25}\ptl_{x_{24}},\eqno(2.76)$$
$$E_{-\al_5}|_V=x_5\ptl_{x_4}+x_7\ptl_{x_6}+x_{13}\ptl_{x_{12}}
+x_{16}\ptl_{x_{15}}
+x_{19}\ptl_{x_{17}}-x_{26}\ptl_{x_{25}},\eqno(2.77)$$
$$E_{-\al_6}|_V=x_8\ptl_{x_5}+x_{10}\ptl_{x_7}+x_{12}\ptl_{x_9}+x_{15}
\ptl_{x_{11}} +x_{17}\ptl_{14}-x_{27}\ptl_{x_{26}},\eqno(2.78)$$
$$E_{-\al_1-\al_3}|_V=-x_3\ptl_{x_1}+x_{14}\ptl_{x_9}+x_{17}\ptl_{x_{12}}+x_{19}\ptl_{x_{13}}
-x_{22}\ptl_{x_{18}}-x_{24}\ptl_{x_{20}},\eqno(2.79)$$
$$E_{-\al_2-\al_4}|_V=x_6\ptl_{x_3}-x_9\ptl_{x_5}-x_{12}\ptl_{x_8}
-x_{20}\ptl_{x_{16}}-x_{23}\ptl_{x_{19}}-x_{25}\ptl_{x_{22}},\eqno(2.80)$$
$$E_{-\al_3-\al_4}|_V=-x_4\ptl_{x_2}-x_{11}\ptl_{x_7}-x_{15}\ptl_{x_{10}}
+x_{18}\ptl_{x_{13}}-x_{22}\ptl_{x_{19}}-x_{25}\ptl_{x_{23}},\eqno(2.81)$$
$$E_{-\al_4-\al_5}|_V=-x_5\ptl_{x_3}+x_9\ptl_{x_6}-x_{13}\ptl_{x_{10}}
+x_{18}\ptl_{x_{15}}+x_{21}\ptl_{x_{17}}-x_{26}\ptl_{x_{24}},\eqno(2.82)$$
$$E_{-\al_5-\al_6}|_V=-x_8\ptl_{x_4}-x_{10}\ptl_{x_6}+x_{13}\ptl_{x_9}+x_{16}
\ptl_{x_{11}} +x_{19}\ptl_{14}-x_{27}\ptl_{x_{25}},\eqno(2.83)$$
$$E_{-\al_1-\al_3-\al_4}|_V=x_4\ptl_{x_1}+x_{14}\ptl_{x_7}+x_{17}\ptl_{x_{10}}
-x_{21}\ptl_{x_{13}}-x_{22}\ptl_{x_{16}}-x_{25}\ptl_{x_{20}},\eqno(2.84)$$
$$E_{-\al_2-\al_3-\al_4}|_V=-x_6\ptl_{x_2}+x_{11}\ptl_{x_5}+x_{15}\ptl_{x_8}
-x_{20}\ptl_{x_{13}}+x_{24}\ptl_{x_{19}}-x_{25}\ptl_{x_{21}},\eqno(2.85)$$
$$E_{-\al_2-\al_4-\al_5}|_V=-x_7\ptl_{x_3}-x_9\ptl_{x_4}+x_{13}\ptl_{x_8}
-x_{20}\ptl_{x_{15}}-x_{23}\ptl_{x_{17}}-x_{26}\ptl_{x_{22}},\eqno(2.86)$$
$$E_{-\al_3-\al_4-\al_5}|_V=x_5\ptl_{x_2}-x_{11}\ptl_{x_6}+x_{16}\ptl_{x_{10}}
+x_{18}\ptl_{x_{12}}-x_{22}\ptl_{x_{17}}-x_{26}\ptl_{x_{23}},\eqno(2.87)$$
$$E_{-\al_4-\al_5-\al_6}|_V=x_8\ptl_{x_3}-x_{12}\ptl_{x_6}-x_{13}\ptl_{x_7}+x_{18}
\ptl_{x_{11}} +x_{21}\ptl_{14}-x_{27}\ptl_{x_{24}},\eqno(2.88)$$
$$E_{-\sum_{i=1}^4\al_i}|_V=x_6\ptl_{x_1}-x_{14}\ptl_{x_5}-x_{17}\ptl_{x_8}
+x_{23}\ptl_{x_{13}}+x_{24}\ptl_{x_{16}}-x_{25}\ptl_{x_{18}},\eqno(2.89)$$
$$E_{-\al_1-\sum_{i=3}^5\al_i}|_V=-x_5\ptl_{x_1}+x_{14}\ptl_{x_6}-x_{19}\ptl_{x_{10}}
-x_{21}\ptl_{x_{12}}-x_{22}\ptl_{x_{15}}-x_{26}\ptl_{x_{20}},\eqno(2.90)$$
$$E_{-\sum_{i=2}^5\al_i}|_V=x_7\ptl_{x_2}+x_{11}\ptl_{x_4}-x_{16}\ptl_{x_8}
-x_{20}\ptl_{x_{12}}+x_{24}\ptl_{x_{17}}-x_{26}\ptl_{x_{21}},\eqno(2.91)$$
$$E_{-\al_2-\sum_{i=4}^6\al_i}|_V=x_{10}\ptl_{x_3}+x_{12}\ptl_{x_4}+x_{13}\ptl_{x_5}-x_{20}
\ptl_{x_{11}}-x_{23}\ptl_{14}-x_{27}\ptl_{x_{22}},\eqno(2.92)$$
$$E_{-\sum_{i=3}^6\al_i}|_V=-x_8\ptl_{x_2}+x_{15}\ptl_{x_6}+x_{16}\ptl_{x_7}+x_{18}
\ptl_{x_9}-x_{22}\ptl_{x_{14}}-x_{27}\ptl_{x_{23}},\eqno(2.93)$$
$$E_{-\sum_{i=1}^5\al_i}|_V=-x_7\ptl_{x_1}-x_{14}\ptl_{x_4}+x_{19}\ptl_{x_8}
+x_{23}\ptl_{x_{12}}+x_{24}\ptl_{x_{15}}-x_{26}\ptl_{x_{18}},\eqno(2.94)$$
$$E_{-\al_1-\sum_{i=3}^6\al_i}|_V=x_8\ptl_{x_1}-x_{17}\ptl_{x_6}-x_{19}\ptl_{x_7}-x_{21}
\ptl_{x_9}-x_{22}\ptl_{x_{11}}-x_{27}\ptl_{x_{20}},\eqno(2.95)$$
$$E_{-\al_4-\sum_{i=2}^5\al_i}|_V=-x_9\ptl_{x_2}+x_{11}\ptl_{x_3}+x_{18}\ptl_{x_8}
-x_{20}\ptl_{x_{10}}+x_{25}\ptl_{x_{17}}-x_{26}\ptl_{x_{19}},\eqno(2.96)$$
$$E_{-\sum_{i=2}^6\al_i}|_V=-x_{10}\ptl_{x_2}-x_{15}\ptl_{x_4}-x_{16}\ptl_{x_5}
-x_{20}\ptl_{x_9}+x_{24}\ptl_{x_{14}}-x_{27}\ptl_{x_{21}},\eqno(2.97)$$
$$E_{-\al_4-\sum_{i=1}^5\al_i}|_V=x_9\ptl_{x_1}-x_{14}\ptl_{x_3}-x_{21}\ptl_{x_8}
+x_{23}\ptl_{x_{10}}+x_{25}\ptl_{x_{15}}-x_{26}\ptl_{x_{16}},\eqno(2.98)$$
$$E_{-\sum_{i=1}^6\al_i}|_V=x_{10}\ptl_{x_1}+x_{17}\ptl_{x_4}+x_{19}\ptl_{x_5}
+x_{23}\ptl_{x_9}+x_{24}\ptl_{x_{11}}-x_{27}\ptl_{x_{18}},\eqno(2.99)$$
$$E_{-\al_4-\sum_{i=2}^6\al_i}|_V=x_{12}\ptl_{x_2}-x_{15}\ptl_{x_3}+x_{18}\ptl_{x_5}
-x_{20}\ptl_{x_7}+x_{25}\ptl_{x_{14}}-x_{27}\ptl_{x_{19}},\eqno(2.100)$$
$$E_{-\al_3-\al_4-\sum_{i=1}^5\al_i}|_V=-x_{11}\ptl_{x_1}+x_{14}\ptl_{x_2}+x_{22}\ptl_{x_8}
-x_{24}\ptl_{x_{10}}+x_{25}\ptl_{x_{12}}-x_{26}\ptl_{x_{13}},\eqno(2.101)$$
$$E_{-\al_4-\sum_{i=1}^6\al_i}|_V=-x_{12}\ptl_{x_1}+x_{17}\ptl_{x_3}-x_{21}\ptl_{x_5}
+x_{23}\ptl_{x_7}+x_{25}\ptl_{x_{11}}-x_{27}\ptl_{x_{16}},\eqno(2.102)$$
$$E_{-\al_4-\al_5-\sum_{i=2}^6\al_i}|_V=-x_{13}\ptl_{x_2}+x_{16}\ptl_{x_3}+x_{18}\ptl_{x_4}
-x_{20}\ptl_{x_6}+x_{26}\ptl_{x_{14}}-x_{27}\ptl_{x_{17}},\eqno(2.103)$$
$$E_{-\al_3-\al_4-\sum_{i=1}^6\al_i}|_V=x_{15}\ptl_{x_1}-x_{17}\ptl_{x_2}+x_{22}\ptl_{x_5}
-x_{24}\ptl_{x_7}+x_{25}\ptl_{x_9}-x_{27}\ptl_{x_{13}},\eqno(2.104)$$
$$E_{-\al_4-\al_5-\sum_{i=1}^6\al_i}|_V=x_{13}\ptl_{x_1}-x_{19}\ptl_{x_3}-x_{21}\ptl_{x_4}
+x_{23}\ptl_{x_6}+x_{26}\ptl_{x_{11}}-x_{27}\ptl_{x_{15}},\eqno(2.105)$$
$$E_{-\sum_{r=3}^5\al_r-\sum_{i=1}^6\al_i}|_V=-x_{16}\ptl_{x_1}+x_{19}\ptl_{x_2}+x_{22}\ptl_{x_4}
-x_{24}\ptl_{x_6}+x_{26}\ptl_{x_9}-x_{27}\ptl_{x_{12}},\eqno(2.106)$$
$$E_{-\al_4-\sum_{r=3}^5\al_r-\sum_{i=1}^6\al_i}|_V=x_{18}\ptl_{x_1}-x_{21}\ptl_{x_2}+x_{22}\ptl_{x_3}
-x_{25}\ptl_{x_6}+x_{26}\ptl_{x_7}-x_{27}\ptl_{x_{10}},\eqno(2.107)$$
$$E_{-\al_2-\al_4-\sum_{r=3}^5\al_r-\sum_{i=1}^6\al_i}|_V=x_{20}\ptl_{x_1}-x_{23}\ptl_{x_2}+x_{24}\ptl_{x_3}
-x_{25}\ptl_{x_4}+x_{26}\ptl_{x_5}-x_{27}\ptl_{x_8}.\eqno(2.108)$$

\section{Proof of the Main Theorem}

Now ${\cal A}=\mbb{C}[x_1,...,x_{27}]$ becomes a ${\cal
G}^{E_6}$-module via the differential operators in (2.35)-(2.108)

According to Table 1, we look for a singular vector of the form:
$$\zeta_1=c_1x_1x_{14}+c_2x_2x_{11}+c_3x_3x_9+c_4x_4x_7+c_5x_5x_6.
\eqno(3.1)$$ By (2.35),
$$0=E_{\al_1}(\zeta_1)=(c_1-c_2)x_1x_{11}\lra c_1=c_2.\eqno(3.2)$$
Moreover, (2.36) implies
$$0=E_{\al_2(\zeta_1)}=-(c_4+c_5)x_4x_5\lra c_5=-c_4.\eqno(3.3)$$
Furthermore, (2.37) gives
$$0=E_{\al_3(\zeta_1)}=(c_2-c_3)x_2x_9\lra c_2=c_3.\eqno(3.4)$$
In addition, (2.38) yields
$$0=E_{\al_4(\zeta_1)}=-(c_3+c_4)x_3x_7\lra c_4=-c_3.\eqno(3.5)$$
The last equation in (3.3) implies $E_{\al_5}(\zeta_1)=0$ by (2.39).
Besides, $E_{\al_6}(\zeta_1)=0$ naturally holds by (3.1) and (2.40).
Taking $c_1=1$, we have the singular vector
$$\zeta_1=x_1x_{14}+x_2x_{11}+x_3x_9-x_4x_7+x_5x_6
\eqno(3.6)$$ of weight $\lmd_6$.

According to (2.73)-(2.78), we set
$$\zeta_2=E_{-\al_6}(\zeta_1)=x_1x_{17}+x_2x_{15}+x_3x_{12}-x_4x_{10}+x_6x_8,
\eqno(3.7)$$
$$\zeta_3=E_{-\al_5}(\zeta_2)=x_1x_{19}+x_2x_{16}+x_3x_{13}-x_5x_{10}
+x_7x_8, \eqno(3.8)$$
$$\zeta_4=E_{-\al_4}(\zeta_3)=x_1x_{21}+x_2x_{18}+x_4x_{13}-x_5x_{12}
+x_8x_9, \eqno(3.9)$$
$$\zeta_5=E_{-\al_3}(\zeta_4)=-x_1x_{22}+x_3x_{18}-x_4x_{16}+x_5x_{15}
-x_8x_{11}, \eqno(3.10)$$
$$\zeta_6=E_{-\al_2}(\zeta_4)=-x_1x_{23}-x_2x_{20}+x_6x_{13}-x_7x_{12}
+x_9x_{10}, \eqno(3.11)$$
$$\zeta_7=E_{-\al_3}(\zeta_6)=x_1x_{24}-x_3x_{20}-x_6x_{16}+x_7x_{15}
-x_{10}x_{11}, \eqno(3.12)$$
$$\zeta_8=E_{-\al_1}(\zeta_5)=-x_2x_{22}-x_3x_{21}+x_4x_{19}-x_5x_{17}
+x_8x_{14}, \eqno(3.13)$$
$$\zeta_9=E_{-\al_4}(\zeta_7)=-x_1x_{25}-x_4x_{20}-x_6x_{18}+x_9x_{15}
-x_{11}x_{12}, \eqno(3.14)$$
$$\zeta_{10}=E_{-\al_1}(\zeta_7)=x_2x_{24}+x_3x_{23}+x_6x_{19}-x_7x_{17}
+x_{10}x_{14}, \eqno(3.15)$$
$$\zeta_{11}=-E_{-\al_5}(\zeta_9)=-x_1x_{26}+x_5x_{20}+x_7x_{18}-x_9x_{16}
+x_{11}x_{13}, \eqno(3.16)$$
$$\zeta_{12}=E_{-\al_4}(\zeta_{10})=-x_2x_{25}+x_4x_{23}+x_6x_{21}-x_9x_{17}
+x_{12}x_{14}, \eqno(3.17)$$
$$\zeta_{13}=E_{-\al_3}(\zeta_{12})=-x_3x_{25}-x_4x_{24}-x_6x_{22}+x_{11}x_{17}
-x_{14}x_{15}, \eqno(3.18)$$
$$\zeta_{14}=-E_{-\al_6}(\zeta_{11})=-x_1x_{27}-x_8x_{20}-x_{10}x_{18}+x_{12}x_{16}
-x_{13}x_{15}, \eqno(3.19)$$
$$\zeta_{15}=-E_{-\al_5}(\zeta_{12})=
-x_2x_{26}-x_5x_{23}-x_7x_{21}+x_9x_{19} -x_{13}x_{14}.
\eqno(3.20)$$

 Define a map $\iota:\ol{1,27}\rta\ol{1,27}$ by
$$\iota(13)=13,\qquad\iota(14)=14,\qquad\iota(15)=15,\eqno(3.21)$$
$$\iota(i)=28-i\qquad\for\;\;i\in\ol{1,27}\setminus\{13,14,15\}.\eqno(3.22)$$
Let $\tau$ be an algebraic automorphism of ${\cal A}$ determined by
$$\tau(x_i)=x_{\iota(i)}\qquad\for\;\;i\in\ol{1,27}.
\eqno(3.23)$$ Now we set
$$\zeta_i=\tau(\zeta_{28-i})\qquad\for\;\;i\in\ol{16,27}.\eqno(3.24)$$
It can be verified that
$$\bar{V}=\sum_{r=1}^{27}\mbb{C}\zeta_r\eqno(3.25)$$ an irreducible ${\cal G}^{E_6}$-submodule and $\{\zeta_r\mid
r\in\ol{1,27}\}$ forms a basis of $\bar V$.

From the Dynkin diagram of $E_6$, we have the following automorphism
of $Q_{E_6}$:
$$\sgm(\sum_{i=1}^6k_i\al_i)=k_6\al_1+k_2\al_2+k_5\al_3+k_4\al_4
+k_3\al_5+k_1\al_6\eqno(3.26)$$ for $\sum_{i=1}^6k_i\al_i\in
Q_{E_6}$. Let $\nu$ be an associative algebra homomorphism of the
associative algebra
$$\mbb{A}=\sum_{i_1,...,i_{27}=0}^\infty {\cal
A}\ptl_{x_1}^{i_1}\cdots\ptl_{x_{27}}^{i_{27}}\eqno(3.27)$$ of
differential operators to itself determined by
$$\nu(x_i)=\zeta_i,\qquad
\nu(\ptl_{x_i})=\ptl_{\zeta_i}\qquad\for\;\;i\in\ol{1,27}.\eqno(3.28)$$
It can be proved that
$$E_{\al}|_{\bar V}=\nu(E_{\sgm(\al)}|_V)\qquad\for\;\;\al\in\Phi_{E_6}^+.
\eqno(3.29)$$ Moreover,
$$\al_j|_{\bar
V}=\sum_{i=1}^{27}b_{i,j}\zeta_i\ptl_{\zeta_i},\eqno(3.30)$$ where
$$b_{i,1}=a_{i,6},\qquad b_{i,3}=a_{i,5},\qquad b_{i,2}=a_{i,2},\qquad b_{i,4}=
a_{i,4}.\eqno(3.31)$$ Thus we have the following table:
\begin{center}{\bf \large Table 2}\end{center}
\begin{center}\begin{tabular}{|r||r|r|r|r|r|r||r||r|r|r|r|r|r|}\hline
$i$&$b_{i,1}$&$b_{i,2}$&$b_{i,3}$&$b_{i,4}$&$b_{i,5}$&$b_{i,6}$&$i$&$b_{i,1}$&$b_{i,2}$&$b_{i,3}$&$b_{i,4}$&$b_{i,5}$&
$b_{i,6}$
\\\hline\hline 1&0&0&0&0&0&1&2&0&0&0&0&1&$-1$\\\hline 3&0&0&0&1&$-1$&0&4&$0$&1&1&$-1$&0&0
\\\hline 5&1&1&$-1$&0&0&0&6&$0$&$-1$&1&$0$&0&0 \\\hline
7&1&$-1$&$-1$&1&0&0&8&$-1$&1&0&$0$&0&0\\\hline
9&1&$0$&0&$-1$&1&0&10&$-1$&$-1$&0&$1$&0&0\\\hline
11&1&$0$&0&$0$&$-1$&1&12&$-1$&0&1&$-1$&1&0\\\hline
13&0&0&$-1$&0&1&0&14&1&0&0&0&0&$-1$\\\hline
15&$-1$&0&1&0&$-1$&1&16&0&0&$-1$&1&$-1$&1\\\hline
17&$-1$&0&1&0&0&$-1$&18&0&1&$0$&$-1$&0&1\\\hline
19&0&0&$-1$&1&0&$-1$&20&0&$-1$&$0$&0&0&1\\\hline
21&0&1&0&$-1$&1&$-1$&22&0&1&0&0&$-1$&0\\\hline
23&0&$-1$&0&0&1&$-1$&24&0&$-1$&0&1&$-1$&0\\\hline
25&0&0&1&$-1$&0&0&26&1&0&$-1$&$0$&0&0\\\hline 27&$-1$&0&0&0&0&0
&&&&&&&\\\hline\end{tabular}\end{center}

According to Table 1 and Table 2, we look for an invariant of the
form
$$\eta=\sum_{i=1}^{12}(d_ix_i\zeta_{28-i}+d_{28-i}x_{28-i}\zeta_i)+d_{13}x_{13}\zeta_{13}
+d_{14}x_{14}\zeta_{14}+d_{15}x_{15}\zeta_{15},\eqno(3.32)$$ where
$d_i\in\mbb{C}$. By (2.35), (2.40) and (3.29), we have
\begin{eqnarray*}& &0=E_{\al_1}(\eta)\\&=&-d_2x_1\zeta_{26}+d_{14}x_{11}\zeta_{14}+d_{17}x_{15}\zeta_{11}+d_{19}x_{16}\zeta_9
+d_{21}x_{18}\zeta_7+d_{23}x_{20}\zeta_5\\ &
&-d_{20}x_{20}\zeta_5-d_{18}x_{18}\zeta_7-d_{16}x_{16}\zeta_9-d_{15}x_{15}\zeta_{11}-d_{11}x_{11}\zeta_{14}
+d_1x_1\zeta_{26},\hspace{2cm}(3.33)
\end{eqnarray*}
\begin{eqnarray*}&
&0=E_{\al_6}(\eta)\\&=&-d_8x_5\zeta_{20}-d_{10}x_7\zeta_{18}-d_{12}x_9\zeta_{16}-d_{15}x_{11}\zeta_{15}
-d_{17}x_{14}\zeta_{11}+d_{27}x_{26}\zeta_1\\ &
&-d_{26}x_{26}\zeta_1+d_{14}x_{14}\zeta_{11}+d_{11}x_{11}\zeta_{15}+d_9x_9\zeta_{16}+d_7x_7\zeta_{18}
+d_5x_5\zeta_{20}.\hspace{2.5cm}(3.34)
\end{eqnarray*}
So we take
$$d_2=d_1,\;\;d_{14}=d_{11},\;\;d_{17}=d_{15},\;\;d_{19}=d_{16},\;\;d_{21}=d_{18},\;\;d_{23}=d_{20},\eqno(3.35)$$
$$d_8=d_7,\;\;d_{10}=d_7,\;\;d_{12}=d_9,\;\;d_{15}=d_{11},\;\;d_{17}=d_{14},\;\;d_{27}=d_{26}.\eqno(3.36)$$
Moreover, (2.37), (2.39) and (3.29) imply
\begin{eqnarray*}& &0=E_{\al_3}(\eta)\\&=&-d_3x_2\zeta_{25}+d_{11}x_9\zeta_{17}+d_{15}x_{12}\zeta_{15}+d_{16}x_{13}
\zeta_{12} +d_{22}x_{21}\zeta_6+d_{24}x_{23}\zeta_4\\ &
&-d_{23}x_{23}\zeta_4-d_{21}x_{21}\zeta_6-d_{13}x_{13}\zeta_{12}-d_{12}x_{12}\zeta_{15}-d_9x_9\zeta_{17}
+d_2x_2\zeta_{25},\hspace{2.2cm}(3.37)
\end{eqnarray*}
\begin{eqnarray*}&
&0=E_{\al_5}(\eta)\\&=&-d_5x_4\zeta_{23}-d_7x_6\zeta_{21}-d_{13}x_{12}\zeta_{13}-d_{16}x_{15}\zeta_{12}
-d_{19}x_{17}\zeta_9+d_{26}x_{25}\zeta_2\\ &
&-d_{25}x_{25}\zeta_2+d_{17}x_{17}\zeta_9+d_{15}x_{15}\zeta_{12}+d_{12}x_{12}\zeta_{13}+d_6x_6\zeta_{21}
+d_4x_4\zeta_{23}.\hspace{2.3cm}(3.38)
\end{eqnarray*}
Hence we get
$$d_3=d_2,\;\;d_{11}=d_9,\;\;d_{15}=d_{12},\;\;d_{16}=d_{13},\;\;d_{22}=d_{21},\;\;d_{24}=d_{23},\eqno(3.39)$$
$$d_5=d_4,\;\;d_7=d_6,\;\;d_{13}=d_{12},\;\;d_{16}=d_{15},\;\;d_{19}=d_{17},
\;\;d_{26}=d_{25}.\eqno(3.40)$$ Furthermore, (2.36), (2.38) and
(3.29) yield
\begin{eqnarray*}&
&0=E_{\al_2}(\eta)\\&=&-d_6x_4\zeta_{22}-d_{22}x_{22}\zeta_4-d_7x_5\zeta_{21}
-d_{21}x_{21}\zeta_5-d_{10}x_8\zeta_{18}-d_{18}x_{18}\zeta_8
\\ &
&+d_{20}x_{18}\zeta_8+d_8x_8\zeta_{18}+d_{23}x_{21}\zeta_5+d_5x_5\zeta_{21}
+d_{24}x_{22}\zeta_4 +d_4x_4\zeta_{22},\hspace{2.7cm}(3.41)
\end{eqnarray*}
\begin{eqnarray*}&
&0=E_{\al_4}(\eta)\\&=&-d_4x_3\zeta_{24}-d_{24}x_{24}\zeta_3-d_9x_7\zeta_{19}
-d_{19}x_{19}\zeta_7-d_{12}x_{10}\zeta_{16}-d_{16}x_{16}\zeta_{10}
\\ &&-d_{18}x_{16}\zeta_{10}-d_{10}x_{10}\zeta_{16}
-d_{21}x_{19}\zeta_7-d_7x_7\zeta_{19} +d_{25}x_{24}\zeta_3
+d_3x_3\zeta_{24}.\hspace{2.3cm}(3.42)
\end{eqnarray*}
Thus we obtain
$$d_6=d_4,\;\;d_{24}=d_{22},\;\;d_7=d_5,\;\;d_{23}=d_{21},\;\;d_{10}=d_8,\;\;
d_{20}=d_{18},\eqno(3.43)$$
$$d_4=d_3,\;\;d_{25}=d_{24},\;\;d_9=-d_7,\;\;d_{21}=-d_{19},\;\;d_{12}=-d_{10},
\;\;d_{18}=-d_{16}.\eqno(3.44)$$ By (3.35), (3.36), (3.39), (3.40),
(3.42) and (3.43), we have
\begin{eqnarray*}\hspace{2cm}&
&d_1=d_2=d_3=d_4=d_5=d_6=d_7=d_8=-d_9=d_{10}=-d_{11}\\
& &=-d_{12}=-d_{13}=-d_{14}
=-d_{15}=-d_{16}=-d_{17}=d_{18}\\
& &=-d_{19}=d_{20}=d_{21}=d_{22}=d_{23}=d_{24}
=d_{25}=d_{26}=d_{27}.\hspace{2.1cm}(3.45)\end{eqnarray*} Therefore,
we have the following invariant
$$\eta=\sum_{i=1}^8(x_i\zeta_{28-i}+x_{28-i}\zeta_i)+x_{10}\zeta_{18}-\sum_{r=9,11,12}
(x_r\zeta_{28-r}+x_{28-r}\zeta_r)-\sum_{s=13}^{15}x_s\zeta_s.\eqno(3.46)$$
According to (3.6)-(3.24),
\begin{eqnarray*}\eta&=&3[(x_1x_{14}+x_2x_{11}+x_3x_9)x_{27}
+(x_1x_{17}+x_2x_{15}+x_3x_{12})x_{26}+(x_1x_{19}+x_2x_{16}\\
&
&+x_3x_{13})x_{25}+(x_4x_{13}-x_5x_{12}+x_8x_9)x_{24}-(x_4x_{16}-x_5x_{15}
+x_8x_{11})x_{23}+(x_6x_{13}\\ & &-x_7x_{12}
+x_9x_{10})x_{22}+(x_7x_{15}+x_6x_{16}
-x_{10}x_{11})x_{21}+(x_4x_{19}-x_5x_{17} +x_8x_{14})x_{20}\\
& &+(x_6x_{18}-x_9x_{15}
+x_{11}x_{12})x_{19}+(x_{10}x_{14}-x_7x_{17})x_{18}+
(x_9x_{16}-x_{11}x_{13})x_{17}\\
& &-x_{12}x_{14}x_{16}+x_{14}x_{15}x_{13}]
+(x_4x_7-x_5x_6)x_{27}+(x_4x_{10}-x_6x_8)x_{26}+(x_5x_{10}
\\ & &-x_7x_8)x_{25}+(x_1x_{21}+x_2x_{18})x_{24}+(x_3x_{18}-x_1x_{22})x_{23}-(x_2x_{22}
+x_3x_{21})x_{20}. \hspace{0.3cm}(3.47)\end{eqnarray*} \pse

{\bf Lemma 3.1}. {\it Any homogeneous singular vector in ${\cal A}$
is a monomial in $x_1,\;\zeta_1$ and $\eta$}.

{\it Proof}. Note that
$$x_1x_{14}=\zeta_1-x_2x_{11}-x_3x_9+x_4x_7-x_5x_6
\eqno(3.48)$$
$$x_1x_{17}=\zeta_2-x_2x_{15}-x_3x_{12}+x_4x_{10}-x_6x_8,
\eqno(3.49)$$
$$x_1x_{19}=\zeta_3-x_2x_{16}-x_3x_{13}+x_5x_{10}
-x_7x_8, \eqno(3.50)$$
$$x_1x_{21}=\zeta_4-x_2x_{18}-x_4x_{13}+x_5x_{12}
-x_8x_9, \eqno(3.51)$$
$$x_1x_{22}=-\zeta_5-x_3x_{18}-x_4x_{16}+x_5x_{15}
-x_8x_{11}, \eqno(3.52)$$
$$x_1x_{23}=-\zeta_6-x_2x_{20}+x_6x_{13}-x_7x_{12}
+x_9x_{10}, \eqno(3.53)$$
$$x_1x_{24}=\zeta_7+x_3x_{20}+x_6x_{16}-x_7x_{15}
+x_{10}x_{11}, \eqno(3.54)$$
$$x_1x_{25}=-\zeta_9-x_4x_{20}-x_6x_{18}+x_9x_{15}
-x_{11}x_{12}, \eqno(3.55)$$
$$x_1x_{26}=-\zeta_{11}+x_5x_{20}+x_7x_{18}-x_9x_{16}
+x_{11}x_{13} \eqno(3.56)$$ by (3.6)-(3.12), (3.14) and (3.16).
Moreover, (3.47) can be written as
\begin{eqnarray*}& &(3x_1x_{14}+3x_2x_{11}+3x_3x_9+x_4x_7-x_5x_6)x_{27}
\\ &=&\eta-3[(x_1x_{17}+x_2x_{15}+x_3x_{12})x_{26}+(x_1x_{19}+x_2x_{16}
+x_3x_{13})x_{25}+(x_4x_{13}\\
& &-x_5x_{12}+x_8x_9)x_{24}-(x_4x_{16}-x_5x_{15}
+x_8x_{11})x_{23}+(x_6x_{13}-x_7x_{12} +x_9x_{10})x_{22}\\
& &+(x_7x_{15}+x_6x_{16} -x_{10}x_{11})x_{21}+(x_4x_{19}-x_5x_{17}
+x_8x_{14})x_{20}+(x_6x_{18}\\ & &-x_9x_{15}
+x_{11}x_{12})x_{19}+(x_{10}x_{14}-x_7x_{17})x_{18}+
(x_9x_{16}-x_{11}x_{13})x_{17}\\ & &
-x_{12}x_{14}x_{16}+x_{14}x_{15}x_{13}]
-(x_4x_{10}-x_6x_8)x_{26}-(x_5x_{10}-x_7x_8)x_{25}\\
&&-(x_1x_{21}+x_2x_{18})x_{24}-(x_3x_{18}-x_1x_{22})x_{23}+(x_2x_{22}
+x_3x_{21})x_{20}. \hspace{2.5cm}(3.57)\end{eqnarray*}

Let $f$ be any homogenous singular vector in ${\cal A}$. According
to the above equations, $f$ can be written as a rational function
$f_1$ in
$$\{x_i,\zeta_r,\eta\mid
i\in\{\ol{1,13},15,16,18,20\};\;r\in\{\ol{1,7},9,11\}\}.\eqno(3.58)$$
By (2.63)-(2.70), (3.28) and (3.29),
$$0=E_{\al_3+\al_4+\sum_{i=1}^5\al_i}(f_1)=x_1\ptl_{x_{11}}(f_1),\eqno(3.59)$$
$$0=E_{\al_4+\sum_{i=1}^6\al_i}(f_1)=x_1\ptl_{x_{12}}(f_1),
\eqno(3.60)$$
$$0=E_{\al_4+\al_5+\sum_{i=2}^6\al_i}(f_1)=x_2\ptl_{x_{13}}(f_1)+\zeta_1\ptl_{\zeta_{11}}(f_1)
, \eqno(3.61)$$
$$0=E_{\al_3+\al_4+\sum_{i=1}^6\al_i}(f_1)=-x_1\ptl_{x_{15}}(f_1),
\eqno(3.62)$$
$$0=E_{\al_4+\al_5+\sum_{i=1}^6\al_i}|_V=-x_1\ptl_{x_{13}}(f_1),
\eqno(3.63)$$
$$0=E_{\sum_{r=3}^5\al_r+\sum_{i=1}^6\al_i}|_V=x_1\ptl_{x_{16}}(f_1),
\eqno(3.64)$$
$$0=E_{\al_4+\sum_{r=3}^5\al_r+\sum_{i=1}^6\al_i}|_V=x_1\ptl_{x_{18}}(f_1),
\eqno(3.65)$$
$$0=E_{\al_2+\al_4+\sum_{r=3}^5\al_r+\sum_{i=1}^6\al_i}(f_1)=x_1\ptl_{x_{20}}(f_1).
\eqno(3.66)$$ So $f_1$ is independent of
$x_{11},x_{12},x_{13},x_{15},x_{16}, x_{18},x_{20}$ and
$\zeta_{11}$, that is, $f_1$ is a rational function in
$$\{x_i,\zeta_r,\eta\mid
i\in\ol{1,10};\;r\in\{\ol{1,7},9\}\}.\eqno(3.67)$$

Next (2.56)-(2.62), (3.28) and (3.29) imply that
$$0=E_{\sum_{i=1}^5\al_i}(f_1)=x_1\ptl_{x_7}(f_1),\eqno(3.68)$$
$$0=E_{\al_1+\sum_{i=3}^6\al_i}(f_1)=-x_1\ptl_{x_8}(f_1), \eqno(3.69)$$
$$0=E_{\al_4+\sum_{i=2}^5\al_i}(f_1)=x_2\ptl_{x_9}(f_1)+\zeta_2\ptl_{\zeta_9}(f_1),\eqno(3.70)$$
$$0=E_{\sum_{i=2}^6\al_i}(f_1)=x_2\ptl_{x_{10}}(f_1)+\zeta_1\ptl_{\zeta_7}(f_1),
\eqno(3.71)$$
$$0=E_{\al_4+\sum_{i=1}^5\al_i}(f_1)=-x_1\ptl_{x_9}(f_1),\eqno(3.72)$$
$$0=E_{\sum_{i=1}^6\al_i}(f_1)=-x_1\ptl_{x_{10}}(f_1).
\eqno(3.73)$$  Hence $f_1$ is independent of
$x_7,x_8,x_9,x_{10},\zeta_7$ and $\zeta_9$, that is, $f_1$ is a
rational function in
$$\{x_i,\zeta_r,\eta\mid
i,r\in\ol{1,6}\}.\eqno(3.74)$$

Now (2.41), (2.45)-(2.50), (2.52), (2.55), (3.28) and (3.29) give
that
$$0=E_{\al_1+\al_3}(f_1)=x_1\ptl_{x_3}(f_1),\eqno(3.75)$$
$$0=E_{\al_5+\al_6}(f_1)=\zeta_1\ptl_{\zeta_3}(f_1),\eqno(3.76)$$
$$0=E_{\al_1+\al_3+\al_4}(f_1)=-x_1\ptl_{x_4}(f_1),\eqno(3.77)$$
$$0=E_{\al_2+\al_3+\al_4}(f_1)=x_2\ptl_{x_6}(f_1),\eqno(3.78)$$
$$0=E_{\al_2+\al_4+\al_5}|_V=\zeta_2\ptl_{\zeta_6}(f_1),\eqno(3.79)$$
$$0=E_{\al_3+\al_4+\al_5}(f_1)=-x_2\ptl_{x_5}(f_1)-\zeta_2\ptl_{\zeta_5}(f_1),
\eqno(3.80)$$
$$0=E_{\al_4+\al_5+\al_6}(f_1)=-\zeta_1\ptl_{\zeta_4}(f_1),\eqno(3.81)$$
$$0=E_{\al_1+\sum_{i=3}^5\al_i}(f_1)=x_1\ptl_{x_5}(f_1).\eqno(3.82)$$
Thus $f_1$ is independent of $\{x_i,\zeta_i\mid i\in\ol{3,6}\}$ ,
that is, $f_1$ is a rational function in
$\{x_1,x_2,\zeta_1,\zeta_2,\eta\}.$ Finally, (2.35), (2.40), (3.28)
and (3.29) yield
$$0=E_{\al_1}(f_1)=-x_1\ptl_{x_2}(f_1),\qquad 0=E_{\al_6}(f_1)=-\zeta_1\ptl_{\zeta_2}(f_1).\eqno(3.83)$$
Therefore, $f_1$ is independent of $x_2$ and $\zeta_2$, that is,
$f=f_1$ is a rational function in $x_1,\;\zeta_1$ and $\eta$. By
(3.48) and (3.57), it must be a polynomial in $x_1,\;\zeta_1$ and
$\eta$. Recall that the weights of $x_1,\;\zeta_1$ and $\eta$ are
$\lmd_1,\;\lmd_6$ and $0$, respectively. The homogeneity of $f$
implies that it must be a monomial in $x_1,\;\zeta_1$ and
$\eta.\qquad\Box$\psp

Let $L(m_1,m_2,m_3)$ be the ${\cal G}^{E_6}$-submodule generated by
$x_1^{m_1}\zeta_1^{m_2}\eta^{m_3}$. Then $L(m_1,m_2,m_3)$ is a
finite-dimensional irreducible submodule of highest weight
$m_1\lmd_1+m_2\lmd_6$. By the Weyl's theorem of completely
reducibility and the above lemma, we have
$${\cal A}=\sum_{m_1,m_2,m_2=0}^\infty L(m_1,m_2,m_3).\eqno(3.84)$$
Recall we denote by $V(\lmd)$ the finite-dimensional irreducible
module of highest weight $\lmd$. The above equation implies
$$\frac{1}{(1-q)^{27}}=\frac{1}{1-q^3}\sum_{m_1,m_2=0}^\infty (\dim
V(m_1\lmd_1+m_2\lmd_6))q^{m_1+2m_2}.\eqno(3.85)$$ Equivalently, we
have:\psp

{\bf Lemma 3.2}. {\it The following dimensional property of
irreducible ${\cal G}^{E_6}$-modules holds:}
$$(1-q)^{26}\sum_{m_1,m_2=0}^\infty (\mbox{\it dim}\:
V(m_1\lmd_1+m_2\lmd_6))q^{m_1+2m_2}=1+q+q^2.\eqno(3.86)$$ \pse

Set
$$W=\sum_{i=1}^{27}\mbb{C}\ptl_{x_i}.\eqno(3.87)$$
Then $W$ isomorphic to the module of linear functions on $V$ via
$\ptl_{x_i}(x_j)=\dlt_{i,j}$. Indeed, the linear map determined by
$\ptl_{x_i}\mapsto \zeta_{\iota(i)}$ (cf. (3.21), (3.22)) is a
${\cal G}^{E_6}$-module isomorphism. We define a linear map
$\Im:{\cal A}\rta \mbb{C}[\ptl_{x_1},...,\ptl_{x_{27}}]$ by
$$\Im(x_1^{\al_1}x_2^{\al_2}\cdots x_{27}^{\al_{27}})=\ptl_{x_1}^{\al_1}\ptl_{x_2}^{\al_2}\cdots
\ptl_{x_{27}}^{\al_{27}}.\eqno(3.89)$$ Set
$${\cal D}=\Im(\eta),\qquad{\cal
D}_1=\sum_{i=1}^{27}x_i\ptl_{x_i},\qquad{\cal
D}_2=\sum_{i=1}^{27}\zeta_i\Im(\zeta_i).\eqno(3.89)$$ Then ${\cal
D},\;{\cal D}_1$ and ${\cal D}_2$ are invariant differential
operators, that is,
$$({\cal D}\xi)|_{\cal A}=(\xi{\cal D})|_{\cal
A},\;\;({\cal D}_r\xi)|_{\cal A}=(\xi{\cal D}_r)|_{\cal
A}\qquad\for\;\;\xi\in{\cal G}^{E_6}.\eqno(3.90)$$

Note that Lemma 3.1 implies
$$V^2=L(2,0,0)+L(0,1,0).\eqno(3.91)$$
Symmetrically,
$$W^2=L'(0,2,0)+L'(1,0,0),\eqno(3.92)$$
where $L'(0,2,0)$ is a module generated by the highest weight vector
$\ptl_{27}^2$ with weight $2\lmd_6$ and $L'(1,0,0)$ is a module
generated by the highest weight vector $\Im(\zeta_{27})$ with weight
$\lmd_1$. Thus the subspace of invariants (the trivial submodule) in
$V^2W^2$ is two-dimensional. The trivial submodule of
$L(0,1,0)L'(1,0,0)$ is $\mbb{C}{\cal D}_2$. In $L(2,0,0)L'(0,2,0)$,
there exists an invariant ${\cal D}_3$ with a term
$x_1^2\ptl_{x_1}^2$. So any invariant in $V^2W^2$ must be in
$\mbb{C}{\cal D}_2+\mbb{C}{\cal D}_3$. In particular, the invariant
differential operator
$$[{\cal D},\eta]={\cal D}\eta-\eta{\cal D}=b_0+b_1{\cal
D}_1+b_2{\cal D}_2+b_3{\cal D}_3\eqno(3.93)$$ for some
$b_s\in\mbb{C}$. According to (3.47), $\eta$ does not contain
$x_1^2$. So $b_3=0$. Moreover, (3.47) also implies $b_0=111$.

According to (3.57), the coefficient of $x_{27}\ptl_{x_{27}}$ in
$[{\cal D}_0,\eta]$ must be $11$, which implies $b_1=11$. Observe
that there exists a unique monomial in $\eta$ containing
$x_1x_{14}$, which is $3x_1x_{14}x_{27}$. Thus the coefficient of
$x_1x_{14}\ptl_{x_1}\ptl_{x_{14}}$ in $[{\cal D},\eta]$ must be $9$,
that is, $b_2=9$. So we have: \psp

{\bf Lemma 3.3}. {\it As operators on ${\cal A}$},
$$[{\cal D},\eta]=111+11{\cal
D}_1+9{\cal D}_2.\eqno(3.94)$$ \pse

Let $m_1$ and $m_2$ be nonnegative integers.
 If ${\cal D}(x_1^{m_1}\zeta_1^{m_2})\neq 0$, then it is also a
singular of degree $m_1+2m_2-3$ with the same weight
$m_1\lmd_1+m_2\lmd_6$. But Lemma 3.1 implies that any singular
vector with weight $m_1\lmd_1+m_2\lmd_6$ must has degree $\geq
m_1+2m_2$. This leads a contradiction. Thus
$${\cal
D}(x_1^{m_1}\zeta_1^{m_2})=0\qquad\for\;\;m_1,m_2\in\mbb{N}.\eqno(3.95)$$
Moreover, (3.90) implies
$${\cal D}(L(m_1,m_2,0))=\{0\}\qquad\for\;\;m_1,m_2\in\mbb{N}.\eqno(3.96)$$
Since ${\cal D}_2(x_1^{m_1}\zeta_1^{m_2})$ is also a singular vector
of degree $m_1+2m_2$ with the same weight $m_1\lmd_1+m_2\lmd_6$, we
have
$${\cal
D}_2(x_1^{m_1}\zeta_1^{m_2})=cx_1^{m_1}\zeta_1^{m_2}\eqno(3.97)$$
for some $c\in\mbb{C}$. Let
$$x_i=0\qquad \for\;\;1,14\neq i\in\ol{1,27}\eqno(3.98)$$
in (3.97) and we get
\begin{eqnarray*}cx_1^{m_1+m_2}x_{14}^{m_2}&=&
\lim_{x_i\rta 0;\;8,10\neq
i\in\ol{2,11}}x_1x_{14}(\ptl_{x_1}\ptl_{x_{14}}+\ptl_{x_2}\ptl_{x_{11}}+\ptl_{{x_3}}\ptl_{x_9}
-\ptl_{x_4}\ptl_{x_7}+\ptl_{x_5}\ptl_{x_6})[x_1^{m_1}\\& &\times
(x_1x_{14}+x_2x_{11}+x_3x_9-x_4x_7+x_5x_6)^{m_2}]\\
&=&m_2(m_1+m_2+4)x_1^{m_1+m_2}x_{14}^{m_2}\hspace{6.6cm}(3.99)
\end{eqnarray*}
by (3.6)-(3.24), that is, $c=m_2(m_1+m_2+4)$. We get:\psp

{\bf Lemma 3.4}. {\it For $m_1,m_2\in\mbb{N}$,}
$${\cal D}_2(x_1^{m_1}\zeta_1^{m_2})=m_2(m_1+m_2+4)x_1^{m_1}\zeta_1^{m_2}.\eqno(3.100)$$
\pse

According to Lemma 3.1,
$$V^4=L(4,0,0)+L(2,1,0)+L(1,0,1).\eqno(3.101)$$
Moreover, $L(1,0,1)=\eta V$. Thus the invariants in $V^4W$ are
$\mbb{C}\eta{\cal D}_1$. Hence
$$[{\cal D}_2,\eta]=c_1\eta+c_2\eta{\cal D}_1\qquad\mbox{for
some}\;\;c_1,c_2\in\mbb{C}.\eqno(3.102)$$ Letting the above equation
act on 1, we have
$${\cal D}_2(\eta)=c_1\eta.\eqno(3.103)$$
By (3.6)-(3.24) and (3.47),
\begin{eqnarray*}& &3c_1x_1x_{14}x_{27}=\lim_{x_i\rta 0;\;14\neq
i\in\ol{2,16}}\eta=\lim_{x_i\rta 0;\;14\neq i\in\ol{2,16}}{\cal
D}_2(\eta)\\
&=&3(x_1x_{14}\ptl_{x_1}x_{14}+x_{14}x_{27}\ptl_{x_{14}}\ptl_{x_{27}}+x_1x_{27}\ptl_{x_1}x_{27})(x_1x_{14}x_{27})
=9x_1x_{14}x_{27},\hspace{1.2cm}(3.104)\end{eqnarray*} So $c_1=3$.
Letting (3.102) act on $x_1$, we have:
$${\cal D}_2(\eta x_1)=(3+c_2)\eta x_1.\eqno(3.105)$$
As (3.104),
\begin{eqnarray*}& &3(3+c_2)x_1^2x_{14}x_{27}=\lim_{x_i\rta 0;\;14\neq
i\in\ol{2,16}}(3+c_2)\eta x_1=\lim_{x_i\rta 0;\;14\neq
i\in\ol{2,16}}{\cal
D}_2(\eta x_1)\\
&=&3(x_1x_{14}\ptl_{x_1}x_{14}+x_{14}x_{27}\ptl_{x_{14}}\ptl_{x_{27}}+x_1x_{27}\ptl_{x_1}x_{27})(x_1^2x_{14}x_{27})
=15x_1^2x_{14}x_{27},\hspace{1cm}(3.106)\end{eqnarray*} Hence
$c_2=2$. We get: \psp

{\bf Lemma 3.5}. {\it As operators on ${\cal A}$,}
$$[{\cal D}_2,\eta]=\eta(3+2{\cal D}_1).\eqno(3.107)$$
\pse

For $m,m_1,m_2\in\mbb{N}$ with $m>0$, we have
\begin{eqnarray*}\hspace{1cm}{\cal D}(\eta^mx_1^{m_1}\zeta_1^{m_2})&=&
[m(111+11m_1+m_2(m_1+m_2+26))\\ & &+ \sum_{s=1}^m
s(33+9(3s+m_1+2m_2))]
\eta^{m-1}x_1^{m_1}\zeta_1^{m_2}\neq0\hspace{1.5cm}(3.108)\end{eqnarray*}
by Lemmas 3.3-3.5. According to (3.84) and (3.108), we have: \psp

{\bf Lemma 3.6}. {\it For any $0\neq f\in{\cal A}$},
$${\cal D}(\eta f)\neq 0.\eqno(3.109)$$
\pse

The above lemma implies that
$$\{f\in{\cal A}\mid{\cal D}(f)\}=\sum_{m_1,m_2=0}^\infty
L(m_1,m_2,0).\eqno(3.110)$$ Recall that ${\cal A}_m$ be the subspace
of homogeneous polynomials of degree $m$ in ${\cal A}$. Denote
$$\Phi_m=\{f\in{\cal A}_m\mid {\cal D}(f)=0\}.\eqno(3.111)$$
In summary, we have the following version of the  main theorem. \psp

{\bf Theorem 3.7}. {\it The set
$\{x^{m_1}_1\zeta_1^{m_2}\eta^{m_3}\mid n_1,m_2,m_3\in\mbb{N}\}$ is
the set of all singular vectors in ${\cal A}$ up to a scalar
multiple. In particular, $\eta$ is the unique fundamental invariant
(up to constant) and  the identity
$$(1-q)^{26}\sum_{m_1,m_2=0}^\infty (\mbox{\it dim}\:
V(m_1\lmd_1+m_2\lmd_6))q^{m_1+2m_2}=1+q+q^2\eqno(3.112)$$ holds.
Furthermore,
$${\cal A}_k=\Phi_k\oplus \eta{\cal A}_{k-3}\qquad\mbox{\it
for}\;\; k\in\mbb{N}\eqno(3.113)$$ and
$$\Phi_m=\sum_{i=0}^{[\!|m/2|\!]}L(m-2i,i,0)\qquad{\it
for}\;\;m\in\mbb{N},\eqno(3.114)$$ where we treat ${\cal A}_r=\{0\}$
if $r<0$.}

\bibliographystyle{amsplain}

\end{document}